\documentclass[12pt]{amsart}
\usepackage{amssymb}
\usepackage{amsfonts}
\usepackage{color,soul}
\newcommand{\nc}{\newcommand}
\nc{\thref}[1]{Theorem~\ref{theo:#1}}
\nc{\selabel}[1]{\label{sect:#1}}
\nc{\seref}[1]{Section~\ref{sect:#1}}
\nc{\lelabel}[1]{\label{lemm:#1}}
\nc{\leref}[1]{Lemma~\ref{lemm:#1}}
\nc{\prlabel}[1]{\label{prop:#1}}
\nc{\prref}[1]{Proposition~\ref{prop:#1}}
\nc{\colabel}[1]{\label{coro:#1}}
\nc{\coref}[1]{Corollary~\ref{coro:#1}}
\nc{\exlabel}[1]{\label{exam:#1}}
\nc{\exref}[1]{Example~\ref{exam:#1}}
\nc{\delabel}[1]{\label{defi:#1}}
\nc{\deref}[1]{Definition~\ref{defi:#1}}
\nc{\eqlabel}[1]{\label{equa:#1}}
\nc{\relabel}[1]{\label{rema:#1}}
\nc{\reref}[1]{Lemma~\ref{rema:#1}}
\providecommand{\operatorname}[1]{\mathrm{#1}\,}
\nc{\Hom}{\operatorname{Hom}} \nc{\Mor}{\operatorname{Mor}}
\nc{\Aut}{\operatorname{Aut}} \nc{\Ann}{\operatorname{Ann}}
\nc{\Ker}{\operatorname{Ker}} \nc{\Trace}{\operatorname{Trace}}
\nc{\Char}{\operatorname{Char}} \nc{\Mod}{\operatorname{Mod}}
\nc{\End}{\operatorname{End}} \nc{\Spec}{\operatorname{Spec}}
\nc{\Span}{\operatorname{Span}} \nc{\sgn}{\operatorname{sgn}}
\nc{\Id}{\operatorname{Id}} \nc{\Com}{\operatorname{Com}}
\nc{\rank}{\operatorname{rank}}
\nc{\Clausen}{\operatorname{Cl}}

\let\:=\colon

\newtheorem{de}{Definition}[section]
\newtheorem{lm}[de]{Lemma}
\newtheorem{pr}[de]{Proposition}
\newtheorem{co}[de]{Corollary}
\newtheorem{re}[de]{Remark}
\newtheorem{res}[de]{Remarks}
\newtheorem{te}[de]{Theorem}
\newtheorem{ex}[de]{Example}
\newtheorem{exs}[de]{Examples}

\def\bex{\begin{ex}}
\def\eex{\end{ex}}
\def\bexs{\begin{exs}}
\def\eexs{\end{exs}}
\def\bl{\begin{lm}}
\def\el{\end{lm}}
\def\bc{\begin{co}}
\def\ec{\end{co}}
\def\bt{\begin{te}}
\def\et{\end{te}}
\def\bpr{\begin{pr}}
\def\epr{\end{pr}}
\def\br{\begin{re}}
\def\er{\end{re}}
\def\brs{\begin{res}}
\def\ers{\end{res}}
\def\bd{\begin{de}}
\def\ed{\end{de}}
\def\be{\begin{equation}}
\def\ee{\end{equation}}
\def\bea{\begin{eqnarray*}}
\def\eea{\end{eqnarray*}}
\def\bp{\begin{proof}}
\def\ep{\end{proof}}

\def\RR{{\mathbb R}}

\def\NN{{\mathbb N}}

\let\:=\colon

\begin{document}

\title[Generalized rational zeta series for $\zeta(2n)$ and $\zeta(2n+1)$]{Generalized rational zeta series for $\zeta(2n)$ and $\zeta(2n+1)$}

\begin{abstract}
In this paper, we find rational zeta series with $\zeta(2n)$ in terms of $\zeta(2k+1)$ and $\beta(2k)$, the Dirichlet beta function. We then develop a certain family of generalized rational zeta series using the generalized Clausen function and use those results to discover a second family of generalized rational zeta series. As a special case of our results from Theorem 3.1, we prove a conjecture given in 2012 by F.M.S. Lima. Later, we use the same analysis but for the digamma function $\psi(x)$ and negapolygammas $\psi^{(-m)}(x)$. With these, we extract the same two families of generalized rational zeta series with $\zeta(2n+1)$ on the numerator rather than $\zeta(2n)$.
\end{abstract}

\author{Derek Orr}



\maketitle

\tableofcontents

\section{Introduction}

\vspace{0.3cm}

In 1734, Leonard Euler proved an amazing result, now known as the celebrated Euler series:

$$\zeta(2) = \sum_{n=1}^{\infty} \frac{1}{n^2} = 1 + \frac{1}{4} + \frac{1}{9} + \frac{1}{16} + \dots = \frac{\pi^2}{6}, $$

\vspace{0.5cm}

where $\zeta(s)$ is the Riemann zeta function, defined as

$$ \zeta(s) = \sum_{n=1}^{\infty} \frac{1}{n^s}, \hspace{5pt} \Re(s) > 1. $$

Later, Euler gave the formula

$$ \zeta(2k)= \sum_{n=1}^{\infty} \frac{1}{n^{2k}} = \frac{(-1)^{k+1}B_{2k}(2\pi)^{2k}}{2(2k)!}, \hspace{5pt} k \in \mathbb{N}_{0},$$

\bigskip

where $B_{n}$ are the Bernoulli numbers, defined by

$$\displaystyle \frac{z}{e^z-1}=\sum_{n=0}^{\infty}\frac{B_{n}}{n!}z^n, \hspace{5pt} |z|<2\pi.$$

\vspace{0.5cm}

These Bernoulli numbers also arise in certain power series, namely of the cotangent function

\begin{equation}
\cot(x)=\sum_{n=0}^{\infty}\frac{(-1)^n2^{2n}B_{2n}}{(2n)!}x^{2n-1} = -2\sum_{n=0}^{\infty} \frac{\zeta(2n)}{\pi^{2n}}x^{2n-1}, \hspace{5pt} |x|<\pi.
\end{equation}

\vspace{0.4cm}

The odd arguments of $\zeta(s)$ are the interesting as they do not have a closed form, though many mathematicians have studied them in detail. In 1979, Roger Apery \cite{Apery} proved that $\zeta(3)$ is irrational using the fast converging series

$$ \zeta(3) = \frac{5}{2}\sum_{n=1}^{\infty} \frac{(-1)^{n-1}}{n^3\binom{2n}{n}}.$$

\vspace{0.4cm}

It is still unknown whether $\zeta(5)$ is irrational however, it is known that at least one of $\zeta(5)$, $\zeta(7)$, $\zeta(9)$, and $\zeta(11)$ is irrational (see \cite{Zudilin}). The rational zeta series in this paper will involve odd arguments of $\zeta(s)$ among other things. Another main function needed is the Clausen function (or Clausen's integral),

\begin{equation}
\Clausen_{2}(\theta) := \sum_{k=1}^{\infty} \frac{\sin(k\theta)}{k^2} = -\int_{0}^{\theta} \log\Big(2\sin\Big(\frac{\phi}{2}\Big)\Big) \hspace{3pt} d\phi.
\end{equation}

\vspace{0.4cm}

The Clausen function also has a power series representation which will be used later in the paper. It is given as

\begin{equation} \frac{\Clausen_{2}(\theta)}{\theta} = 1-\log |\theta|+\sum_{n=1}^{\infty}\frac{\zeta(2n)}{n(2n+1)}\Big(\frac{\theta}{2\pi}\Big)^{2n}, \hspace{5pt} |\theta|<2\pi.
\end{equation}

\vspace{0.3cm}

There are also higher order Clausen-type function defined as

\begin{equation}
\Clausen_{2m}(\theta) := \sum_{k=1}^{\infty} \frac{\sin(k\theta)}{k^{2m}}, \hspace{0.5cm} \Clausen_{2m+1}(\theta) := \sum_{k=1}^{\infty} \frac{\cos(k\theta)}{k^{2m+1}}.
\end{equation}

\vspace{0.4cm}

The Clausen function is widely studied and has many applications in mathematics and mathematical physics (\cite{Choi}, \cite{ChoiSriv1}, \cite{Doeder}, \cite{Grosjean}, \cite{OrrLupu}, \cite{SrivGlasAdam}, \cite{WuZhangLiu}). We will also discuss the Dirichlet beta function

$$ \beta(s) = \sum_{n=0}^{\infty} \frac{(-1)^n}{(2n+1)^s}, \hspace{5pt} \Re(s) > 0. $$

\vspace{0.4cm}

When $s=2$,

$$ \beta(2) = G = \sum_{n=0}^{\infty} \frac{(-1)^n}{(2n+1)^2} $$

\vspace{0.4cm}

is known as Catalan's constant. Using this and the Riemann zeta function, we find

\begin{equation}
\Clausen_{2m}(\pi) = 0, \hspace{0.5cm} \Clausen_{2m+1}(\pi) = -\frac{(4^m-1)\zeta(2m+1)}{4^m},
\end{equation}

\vspace{0.4cm}

and

\begin{equation}
\Clausen_{2m}(\pi/2) = \beta(2m), \hspace{0.5cm} \Clausen_{2m+1}(\pi/2) = -\frac{(4^m-1)\zeta(2m+1)}{2^{4m+1}}.
\end{equation}

\vspace{0.5cm}

Also, from $(4)$ we can see

\begin{equation} 
\frac{d}{d\theta} \Clausen_{2m}(\theta) = \Clausen_{2m-1}(\theta), \hspace{0.5cm} \frac{d}{d\theta} \Clausen_{2m+1}(\theta) = -\Clausen_{2m}(\theta),
\end{equation}

\begin{equation}
\int_{0}^{\theta} \Clausen_{2m}(x) \hspace{3pt} dx = \zeta(2m+1) - \Clausen_{2m+1}(\theta), \hspace{0.5cm} \int_{0}^{\theta} \Clausen_{2m-1}(x) \hspace{3pt} dx = \Clausen_{2m}(\theta).
\end{equation}

\vspace{0.4cm}

Using $(4)$ and $(7)$, we find

\begin{equation}
\Clausen_{1}(\theta) = -\log\Big(2\sin\Big(\frac{\theta}{2}\Big)\Big), \hspace{5pt} |\theta| < 2\pi.
\end{equation}

\vspace{0.3cm}

Writing others out,

$$\Clausen_{3}(z) = \zeta(3) - \int_{0}^{z} \Clausen_{2}(t) \hspace{3pt} dt,$$

$$\Clausen_{4}(z) = \int_{0}^{z} \Clausen_{3}(x) \hspace{3pt} dx = z\zeta(3)-\int_{0}^{z} \int_{0}^{x} \Clausen_{2}(t) \hspace{3pt} dt \hspace{3pt} dx = z\zeta(3)-\int_{0}^{z} (z-t)\Clausen_{2}(t) \hspace{3pt} dt, $$

\begin{multline*}
\Clausen_{5}(z) = \zeta(5) - \int_{0}^{z} \Clausen_{4}(x) \hspace{3pt} dx = \zeta(5) - \frac{1}{2}z^2\zeta(3)+\frac{1}{2}\int_{0}^{z} (z-t)^2\Clausen_{2}(t) \hspace{3pt} dt,
\end{multline*}

\vspace{0.4cm}

and by induction, for $m \geq 3$,

\begin{multline}
\Clausen_{m}(z) = (-1)^{\lfloor \frac{m-1}{2} \rfloor}\sum_{k=1}^{\lfloor \frac{m-1}{2} \rfloor} \frac{(-1)^kz^{m-2k-1}}{(m-2k-1)!}\zeta(2k+1)\\ + \frac{(-1)^{\lfloor \frac{m-1}{2} \rfloor}}{(m-3)!}\int_{0}^{z}(z-t)^{m-3}\Clausen_{2}(t) \hspace{3pt} dt.
\end{multline}

\vspace{0.4cm}

In the latter sections of this paper, we will focus on the polygamma function. Its definition is given by

$$ \psi^{(n)}(z) := \frac{d^{n+1}}{dz^{n+1}} \log\Gamma(z), \hspace{5pt} n \in \NN_{0}.$$

\vspace{0.4cm}

This paper will discuss when $n=0$. For $n=0$, $\psi^{(0)}(z) = \psi(z)$ is called the digamma function. It has been shown (see \cite{AdamchikPG}) that there is a closed form of $\displaystyle \int_{0}^{z} x^n\psi(x) \hspace{3pt} dx$ in terms of sums involving Harmonic numbers, Bernoulli numbers and Bernoulli polynomials. There are also definitions for negative order polygamma functions, called negapolygammas, given by

$$ \psi^{(-1)}(z) = \log\Gamma(z),$$

$$ \psi^{(-2)}(z) = \int_{0}^{z} \log\Gamma(x) \hspace{3pt} dx,$$

$$ \psi^{(-3)}(z) = \int_{0}^{z}\int_{0}^{x} \log\Gamma(t) \hspace{3pt} dt \hspace{2pt} dx = \int_{0}^{z}\int_{t}^{z} \log\Gamma(t) \hspace{3pt} dx \hspace{2pt} dt = \int_{0}^{z} (z-t)\log\Gamma(t) \hspace{3pt} dt,$$

\vspace{0.3cm}

and by induction, for $n \geq 2$,

\begin{equation}
\psi^{(-n)}(z) = \frac{1}{(n-2)!}\int_{0}^{z} (z-t)^{n-2}\log\Gamma(t) \hspace{3pt} dt.
\end{equation}

\vspace{0.4cm}

Note the Taylor series for $\log\Gamma(z)$ is

\begin{equation}
\log\Gamma(z) = -\log z -\gamma z + \sum_{k=2}^{\infty} \frac{(-1)^k\zeta(k)}{k}z^k, \hspace{5pt} |z| < 1,
\end{equation}

\vspace{0.4cm}

where $\gamma$ is the Euler-Mascheroni constant. The last function to introduce is the Hurwitz zeta function, defined by

$$ \zeta(s,a) = \sum_{k=0}^{\infty} \frac{1}{(k+a)^s}, \hspace{5pt} a \in \RR \backslash -\NN, \hspace{5pt} \Re(s) > 1. $$

\vspace{0.4cm}

The negapolygammas are related to the derivative of the Hurwitz zeta function with respect to the first variable (see \cite{AdamchikPG}).

\subsection{Organization of the Paper}
We begin by investigating $\int_{0}^{\pi z} x^p\cot(x) \hspace{3pt} dx$ for $|z|<1$, which is studied in \cite{SrivGlasAdam}. We then compute this same integral using the power series for the cotangent and obtain a rational zeta series representation with $\zeta(2n)$ on the numerator for $z=1/2$ and $z=1/4$. Afterwords, we work on a generalized rational zeta series with $\zeta(2n)$ on the numerator and an arbitrary number of monomials on the denominator using the generalized Clausen function $\Clausen_{m}(z)$. As a special case, we immediately prove a conjecture given in 2012 by F.M.S. Lima. After doing this, we go back to the cotangent function and discover a separate class of generalized $\zeta(2n)$ series. Lastly, we perform the same analysis but for the digamma function $\psi(x)$ instead of $\cot(x)$ and instead of $\Clausen_{m}(x)$, we use the negapolygammas $\psi^{(-m)}(x)$. Using the $\zeta(2n)$ sums from earlier, we extract the same rational zeta series representations with $\zeta(2n+1)$ on the numerator for $z=1/2$ and $z=1/4$.

\vspace{0.5cm}

\textbf{Acknowledgements.} I would like to thank Cezar Lupu, Tom Hales, and George Sparling for their interest in my paper as well as their advice which led to some improvements of the paper. 

\bigskip

\section{Rational $\zeta(2n)$ series with $\cot(x)$}

\vspace{0.3cm}

\bt{For $p \in \mathbb{N}$ and $|z| < 1$,}

\begin{multline}
\int_{0}^{\pi z} x^p\cot(x) \hspace{3pt} dx = (\pi z)^p\sum_{k=0}^{p} \frac{p!(-1)^{\lfloor \frac{k+3}{2} \rfloor}}{(p-k)!(2\pi z)^k}\Clausen_{k+1}(2\pi z) + \delta_{\lfloor \frac{p}{2} \rfloor, \frac{p}{2}}\frac{p!(-1)^{\frac{p}{2}}}{2^p}\zeta(p+1),
\end{multline}

\vspace{0.3cm}

where $\delta_{j,k}$ is the Kronecker delta function.

\et

\vspace{0.3cm}

\textit{Proof.} Let $f(z)$ be the left hand side of the equation and let $g(z)$ be the right hand side. Note that $f'(z) = \pi^{p+1}z^p\cot(\pi z)$. Using $(7)$ and $(9)$,

\begin{multline*}
g'(z) = \frac{1}{2^p}\sum_{k=0}^{p} \frac{p!(-1)^{\lfloor \frac{k+3}{2} \rfloor}(p-k)(2\pi)^{p-k}z^{p-k-1}}{(p-k)!}\Clausen_{k+1}(2\pi z)+\frac{(\pi z)^p2\pi\cos(\pi z)}{2\sin(\pi z)}\\+\frac{1}{2^p}\sum_{k=1}^{p} \frac{p!(-1)^{\lfloor \frac{k+3}{2} \rfloor}(2\pi z)^{p-k}}{(p-k)!}\Big\{(-1)^{k+1}2\pi\Clausen_{k}(2\pi z)\Big\}
\end{multline*}
$$ = \pi^{p+1}z^p\cot(\pi z) + (\pi z)^p\sum_{k=1}^{p} \frac{p!\Clausen_{k}(2\pi z)}{(p-k)!z^{k+1}(2\pi)^{k-1}}\Big((-1)^{\lfloor \frac{k+2}{2} \rfloor} + (-1)^{\lfloor \frac{k+3}{2} \rfloor}(-1)^k\Big).$$

\vspace{0.4cm}

So indeed, $g'(z) = \pi^{p+1}z^p\cot(\pi z) = f'(z)$. Clearly, $f(0)=0$. For $g(z)$, note that all terms in the sum are zero except when $k=p$. So we have

$$ g(0) = \frac{1}{2^p}p!(-1)^{\lfloor \frac{p+3}{2} \rfloor}\Clausen_{p+1}(0) + \delta_{\lfloor \frac{p}{2} \rfloor, \frac{p}{2}}\frac{p!(-1)^{\frac{p}{2}}}{2^p}\zeta(p+1). $$

\vspace{0.4cm}

From $(4)$, $\Clausen_{p+1}(0) = \delta_{\lfloor \frac{p}{2} \rfloor, \frac{p}{2}}\zeta(p+1)$. So we see $g(0)=0$. Since $f(0)=g(0)$ and $f'(z)=g'(z)$, $f(z)=g(z)$. $\square$

\bigskip

Using $(5)$, $(6)$, and $(9)$, setting $z=1/2$ and $z=1/4$, we find

\begin{multline}
\displaystyle \int_{0}^{\pi/2} x^p\cot(x) \hspace{3pt} dx  = \Big(\frac{\pi}{2}\Big)^p\bigg(\log2+\sum_{k=1}^{\lfloor\frac{p}{2}\rfloor} \frac{p!(-1)^k(4^k-1)}{(p-2k)!(2\pi)^{2k}}\zeta(2k+1)\bigg) + \delta_{\lfloor \frac{p}{2} \rfloor, \frac{p}{2}}\frac{p!(-1)^{\frac{p}{2}}\zeta(p+1)}{2^p},
\end{multline}

\vspace{0.3cm}

and

\begin{multline}
\int_{0}^{\pi/4} x^p\cot(x) \hspace{3pt} dx = \frac{1}{2}\Big(\frac{\pi}{4}\Big)^p\bigg(\log2+\sum_{k=1}^{\lfloor \frac{p}{2} \rfloor} \frac{p!(-1)^k(4^k-1)}{(p-2k)!(2\pi)^{2k}}\zeta(2k+1)\\-\sum_{k=1}^{\lfloor \frac{p+1}{2} \rfloor} \frac{p!(-4)^k\beta(2k)}{(p+1-2k)!\pi^{2k-1}}\bigg) + \delta_{\lfloor \frac{p}{2} \rfloor, \frac{p}{2}}\frac{p!(-1)^{\frac{p}{2}}}{2^p}\zeta(p+1).
\end{multline}

\vspace{0.5cm}

We can also integrate $x^p\cot(x)$ using $(1)$ and Fubini's theorem. Doing so, we obtain

$$ \int_{0}^{\pi z} x^p\cot(x) \hspace{3pt} dx = -2\int_{0}^{\pi z} \sum_{n=0}^{\infty} \frac{\zeta(2n)}{\pi^{2n}}x^{2n-1+p} \hspace{3pt} dx = -2(\pi z)^p\sum_{n=0}^{\infty} \frac{\zeta(2n)z^{2n}}{2n+p}, $$

\vspace{0.3cm}

and using $(13)$,

\begin{equation}
-2\sum_{n=0}^{\infty} \frac{\zeta(2n)z^{2n}}{2n+p} = \sum_{k=0}^{p} \frac{p!(-1)^{\lfloor \frac{k+3}{2} \rfloor}}{(p-k)!(2\pi z)^k}\Clausen_{k+1}(2\pi z) + \delta_{\lfloor \frac{p}{2} \rfloor, \frac{p}{2}}\frac{p!(-1)^{\frac{p}{2}}}{(2\pi z)^p}\zeta(p+1).
\end{equation}

\vspace{0.4cm}

So for $z=1/2$ and $z=1/4$, we have

\begin{equation}
- 2\sum_{n=0}^{\infty} \frac{\zeta(2n)}{(2n+p)4^n} = \log2+\sum_{k=1}^{\lfloor\frac{p}{2}\rfloor} \frac{p!(-1)^k(4^k-1)\zeta(2k+1)}{(p-2k)!(2\pi)^{2k}}+\delta_{\lfloor \frac{p}{2} \rfloor, \frac{p}{2}}\frac{p!(-1)^{\frac{p}{2}}\zeta(p+1)}{\pi^p},
\end{equation}

\vspace{0.3cm}

and

\begin{multline}
- 2\sum_{n=0}^{\infty} \frac{\zeta(2n)}{(2n+p)16^n} = \frac{1}{2}\log2
+\frac{1}{2}\sum_{k=1}^{\lfloor\frac{p}{2}\rfloor} \frac{p!(-1)^k(4^k-1)\zeta(2k+1)}{(p-2k)!(2\pi)^{2k}} \\- \frac{\pi}{2} \sum_{k=1}^{\lfloor\frac{p+1}{2}\rfloor} \frac{p!(-4)^k\beta(2k)}{(p+1-2k)!\pi^{2k}}+\delta_{\lfloor \frac{p}{2} \rfloor, \frac{p}{2}}\frac{p!(-1)^{\frac{p}{2}}2^p\zeta(p+1)}{\pi^p}.
\end{multline}

\vspace{0.4cm}

Plugging in $p=1$ into both of the equations yields two nice series representations for $\log2$,

\begin{equation}
-2\sum_{n=0}^{\infty} \frac{\zeta(2n)}{(2n+1)4^n} = \log2,
\end{equation}

\vspace{0.3cm}

and

\begin{equation}
-\frac{4G}{\pi}-4\sum_{n=0}^{\infty} \frac{\zeta(2n)}{(2n+1)16^n} = \log2.
\end{equation}

\vspace{0.4cm}

Other series from $(17)$ and $(18)$ are

$$ \sum_{n=0}^{\infty} \frac{\zeta(2n)}{(n+1)4^n} = -\log2+\frac{7\zeta(3)}{2\pi^2},$$
$$ \sum_{n=0}^{\infty} \frac{\zeta(2n)}{(2n+3)4^n} = -\frac{1}{2}\log2+\frac{9\zeta(3)}{4\pi^2}, $$
$$ \sum_{n=0}^{\infty} \frac{\zeta(2n)}{(n+2)4^n} = -\log2 + \frac{9\zeta(3)}{\pi^2}-\frac{93\zeta(5)}{2\pi^4}, $$
$$\sum_{n=0}^{\infty} \frac{\zeta(2n)}{(n+1)16^n} = -\frac{1}{2}\log2+\frac{35\zeta(3)}{4\pi^2}-\frac{4G}{\pi}, $$

\vspace{0.3cm}

and

$$ \sum_{n=0}^{\infty} \frac{\zeta(2n)}{(2n+3)16^n}  = -\frac{1}{4}\log2+\frac{9\zeta(3)}{8\pi^2}-\frac{3G}{\pi}+\frac{24\beta(4)}{\pi^3}. $$

\vspace{0.5cm}

\section{General $\zeta(2n)$ series using $\Clausen_{m}(x)$}

\vspace{0.3cm}

\bt{For $p \in \mathbb{N}_{0}$, $m \in \mathbb{N}$ and $|z| < 1$,}

\begin{multline}
\int_{0}^{2\pi z} x^p\Clausen_{m}(x) \hspace{3pt} dx = - \sum_{k=m+1}^{m+p+1} \frac{(2\pi z)^{p+m+1-k}p!(-1)^{\lfloor \frac{m}{2} \rfloor}(-1)^{\lfloor \frac{k}{2} \rfloor}}{(p+m+1-k)!}\Clausen_{k}(2\pi z) \\+ \delta_{\lfloor \frac{p+m}{2} \rfloor, \frac{p+m}{2}}(-1)^{\lfloor \frac{m}{2} \rfloor}p!(-1)^{\frac{p+m}{2}}\zeta(p+m+1).
\end{multline}
\et

\vspace{0.3cm}

\textit{Proof.} Similar to the proof of $(13)$, we will call the left hand side $f(z)$ and the right hand side $g(z)$. Note $f(0)=g(0)$ and $f'(z)=(2\pi)^{p+1}z^p\Clausen_{m}(2\pi z)$. Using $(7)$,

\begin{multline*}
g'(z)= -\sum_{k=m+1}^{m+p} \frac{p!(2\pi)^{p+m+1-k}z^{p+m-k}(-1)^{\lfloor \frac{m}{2} \rfloor + \lfloor \frac{k}{2} \rfloor}}{(p+m-k)!}\Clausen_{k}(2\pi z)\\ - \sum_{k=m+1}^{m+p+1} \frac{p!(2\pi z)^{p+m+1-k}(-1)^{\lfloor \frac{m}{2} \rfloor+\lfloor \frac{k}{2} \rfloor}}{(p+m+1-k)!}\Big\{(-1)^k2\pi\Clausen_{k-1}(2\pi z)\Big\}
\end{multline*}
\begin{multline*}
= -(2\pi)^{p+1}z^p(-1)^{\lfloor \frac{m}{2} \rfloor+\lfloor \frac{m+1}{2} \rfloor}(-1)^{m+1}\Clausen_{m}(2\pi z) - 2\pi\sum_{k=m+1}^{m+p} \frac{p!(2\pi z)^{p+m-k}(-1)^{\lfloor \frac{m}{2} \rfloor}}{(p+m-k)!}*\\\Clausen_{k}(2\pi z)\Big\{(-1)^{\lfloor \frac{k+1}{2} \rfloor}(-1)^{k+1}+(-1)^{\lfloor \frac{k}{2} \rfloor}\Big\} = (2\pi)^{p+1}z^p\Clausen_{m}(2\pi z).
\end{multline*}

\vspace{0.4cm}

So indeed $f'(z)=g'(z)$ and $f(0)=g(0)$, which implies $f(z)=g(z)$. $\square$

\bigskip

We can see if $p=0$, we recover $(8)$. Now letting $z=1/2$ and $z=1/4$, we find

\begin{multline}
\int_{0}^{\pi} x^p\Clausen_{m}(x) \hspace{3pt} dx = (-1)^{\lfloor \frac{m}{2} \rfloor}p!\pi^{p+m}\sum_{k=\lfloor \frac{m+1}{2} \rfloor}^{\lfloor\frac{p+m}{2}\rfloor} \frac{(-1)^k(4^k-1)\zeta(2k+1)}{(p+m-2k)!(2\pi)^{2k}} \\+ \delta_{\lfloor \frac{p+m}{2} \rfloor, \frac{p+m}{2}}(-1)^{\lfloor \frac{m}{2} \rfloor}p!(-1)^{\frac{p+m}{2}}\zeta(p+m+1),
\end{multline}

\vspace{0.3cm}

and

\begin{multline}
\int_{0}^{\pi/2} x^p\Clausen_{m}(x) \hspace{3pt} dx = (-1)^{\lfloor \frac{m}{2} \rfloor}p!\Big(\frac{\pi}{2}\Big)^{p+m}\bigg(\frac{1}{2}\sum_{k=\lfloor \frac{m+1}{2} \rfloor}^{\lfloor\frac{p+m}{2}\rfloor} \frac{(-1)^k(4^k-1)\zeta(2k+1)}{(p+m-2k)!(2\pi)^{2k}}\\- \frac{\pi}{2}\sum_{k=\lfloor \frac{m+2}{2} \rfloor}^{\lfloor\frac{p+m+1}{2}\rfloor} \frac{(-1)^k4^k\beta(2k)}{(p+m+1-2k)!\pi^{2k}}\bigg) + \delta_{\lfloor \frac{p+m}{2} \rfloor, \frac{p+m}{2}}(-1)^{\lfloor \frac{m}{2} \rfloor}p!(-1)^{\frac{p+m}{2}}\zeta(p+m+1).
\end{multline}

\vspace{0.5cm}

Now, we will integrate the left hand side of $(21)$ using $(3)$ and $(10)$, giving us

\begin{multline*}
\int_{0}^{2\pi z} x^p\Clausen_{m}(x) \hspace{3pt} dx = \int_{0}^{2\pi z} x^p\bigg(\sum_{k=1}^{\lfloor \frac{m-1}{2} \rfloor} \frac{(-1)^{\lfloor \frac{m-1}{2} \rfloor + k}x^{m-2k-1}}{(m-2k-1)!}\zeta(2k+1) \\+ \frac{(-1)^{\lfloor \frac{m-1}{2} \rfloor}}{(m-3)!}\int_{0}^{x}(x-t)^{m-3}\Clausen_{2}(t) \hspace{3pt} dt\bigg) \hspace{3pt} dx
\end{multline*}
$$ = \sum_{k=1}^{\lfloor \frac{m-1}{2} \rfloor} \frac{(-1)^{\lfloor \frac{m-1}{2} \rfloor + k}(2\pi z)^{m+p-2k}}{(m-1-2k)!(m+p-2k)}\zeta(2k+1)+\frac{(-1)^{\lfloor \frac{m-1}{2} \rfloor}}{(m-3)!}\int_{0}^{2\pi z} \int_{0}^{x} x^p(x-t)^{m-3}\Clausen_{2}(t) \hspace{3pt} dt \hspace{2pt} dx$$
\begin{multline*}
=  \sum_{k=1}^{\lfloor \frac{m-1}{2} \rfloor} \frac{(-1)^{\lfloor \frac{m-1}{2} \rfloor + k}(2\pi z)^{m+p-2k}}{(m-1-2k)!(m+p-2k)}\zeta(2k+1)+\frac{(-1)^{\lfloor \frac{m-1}{2} \rfloor}}{(m-3)!}*\\\int_{0}^{2\pi z}\int_{0}^{x} x^p(x-t)^{m-3}\bigg(t-t\log t+\sum_{n=1}^{\infty} \frac{\zeta(2n)t^{2n+1}}{n(2n+1)(2\pi)^{2n}}\bigg) \hspace{3pt} dt \hspace{2pt} dx
\end{multline*}
\begin{multline*}
= \sum_{k=1}^{\lfloor \frac{m-1}{2} \rfloor} \frac{(-1)^{\lfloor \frac{m-1}{2} \rfloor + k}(2\pi z)^{m+p-2k}}{(m-1-2k)!(m+p-2k)}\zeta(2k+1)+\frac{(-1)^{\lfloor \frac{m-1}{2} \rfloor}}{(m-3)!}\int_{0}^{2\pi z} x^p\bigg(\frac{x^{m-1}(H_{m-1}-\log x)}{(m-1)(m-2)}\\+\sum_{n=1}^{\infty} \frac{2\zeta(2n)\Gamma(m-2)x^{m+2n-1}}{2n(2n+1)(2n+2)\dots(2n+m-1)(2\pi)^{2n}}\bigg) \hspace{3pt} dx,
\end{multline*}

\vspace{0.3cm}

where $H_{k}$ is the $k$-th harmonic number and $H_{0} := 0$. Integrating again and simplifying, we arrive at

\begin{multline}
\int_{0}^{2\pi z} x^p\Clausen_{m}(x) \hspace{3pt} dx = \frac{(2\pi z)^{p+m}(-1)^{\lfloor \frac{m-1}{2} \rfloor}}{(m-1)!}\bigg(\frac{H_{m-1}-\log(2\pi z)}{(p+m)}+\frac{1}{(p+m)^2}\\+\sum_{k=1}^{\lfloor \frac{m-1}{2} \rfloor} \frac{(m-1)!(-1)^k\zeta(2k+1)}{(m-1-2k)!(p+m-2k)(2\pi z)^{2k}}+\sum_{n=1}^{\infty} \frac{(m-1)!\zeta(2n)z^{2n}}{n(2n+1)\dots(2n+m-1)(2n+p+m)}\bigg).
\end{multline}

\vspace{0.5cm}

Using $(21)$, we can rearrange this and find

\begin{multline}
\sum_{n=1}^{\infty} \frac{\zeta(2n)z^{2n}}{n(2n+1)\dots(2n+m-1)(2n+m+p)} = \frac{\log(2\pi z)-H_{m-1}}{(m-1)!(p+m)} \\+\sum_{k=m}^{m+p} \frac{(-1)^mp!(-1)^{\lfloor \frac{k+1}{2} \rfloor}}{(p+m-k)!(2\pi z)^k}\Clausen_{k+1}(2\pi z)+\delta_{\lfloor \frac{p+m}{2} \rfloor, \frac{p+m}{2}}\frac{(-1)^{m+1}p!(-1)^{\frac{p+m}{2}}}{(2\pi z)^{p+m}}\zeta(p+m+1)\\- \frac{1}{(m-1)!(p+m)^2} - \sum_{k=1}^{\lfloor \frac{m-1}{2} \rfloor} \frac{(-1)^k\zeta(2k+1)}{(m-1-2k)!(m+p-2k)(2\pi z)^{2k}}.
\end{multline}

\vspace{0.3cm}

Note when $p=0$, we have the very special formula

\begin{multline}
\sum_{n=1}^{\infty} \frac{\zeta(2n)z^{2n}}{n(2n+1)\dots(2n+m)} = \frac{(-1)^m(-1)^{\lfloor \frac{m+1}{2} \rfloor}}{(2\pi z)^m}\Clausen_{m+1}(2\pi z)\\-\delta_{\lfloor \frac{m}{2} \rfloor, \frac{m}{2}}\frac{(-1)^{\frac{3m}{2}}}{(2\pi z)^{m}}\zeta(m+1) - \sum_{k=1}^{\lfloor \frac{m-1}{2} \rfloor} \frac{(-1)^k\zeta(2k+1)}{(m-2k)!(2\pi z)^{2k}} + \frac{\log(2\pi z)-H_{m}}{m!}.
\end{multline}

\vspace{0.4cm}

Letting $z=1/2$ and $z=1/4$, we see that

\begin{multline}
\sum_{n=1}^{\infty} \frac{\zeta(2n)}{n(2n+1)\dots(2n+m)4^n} = \delta_{\lfloor \frac{m}{2} \rfloor, \frac{m}{2}}\frac{(-1)^{\frac{3m+2}{2}}(2^{m+1}-1)\zeta(m+1)}{(2\pi)^m}\\-\sum_{k=1}^{\lfloor \frac{m-1}{2} \rfloor} \frac{(-1)^k\zeta(2k+1)}{(m-2k)!\pi^{2k}}+\frac{\log\pi-H_{m}}{m!},
\end{multline}

\vspace{0.3cm}

and 

\begin{multline}
\sum_{n=1}^{\infty} \frac{\zeta(2n)}{n(2n+1)\dots(2n+m)16^n} = \delta_{\lfloor \frac{m+1}{2} \rfloor, \frac{m+1}{2}}\frac{(-1)^{\frac{3m+1}{2}}2^m\beta(m+1)}{\pi^m}\\ -\frac{1}{2}\delta_{\lfloor \frac{m}{2} \rfloor, \frac{m}{2}}\frac{(-1)^{\frac{3m}{2}}(2^{2m+1}+2^m-1)\zeta(m+1)}{(2\pi)^m} -\sum_{k=1}^{\lfloor \frac{m-1}{2} \rfloor} \frac{(-4)^k\zeta(2k+1)}{(m-2k)!\pi^{2k}}+ \frac{\log(\pi/2)-H_{m}}{m!},
\end{multline}

\vspace{0.4cm}

both of which were conjectured in 2012 and only for $m$ odd (see \cite{Lima}). For general $p$, we have

\begin{multline}
\sum_{n=1}^{\infty} \frac{\zeta(2n)}{n(2n+1)\dots(2n+m-1)(2n+m+p)4^n} = \frac{\log\pi-H_{m-1}}{(m-1)!(p+m)}\\+(-1)^{m+1}p!\bigg(\delta_{\lfloor \frac{p+m}{2} \rfloor, \frac{p+m}{2}}\frac{(-1)^{\frac{p+m}{2}}}{\pi^{p+m}}\zeta(p+m+1) + \sum_{k=\lfloor \frac{m+1}{2} \rfloor}^{\lfloor\frac{p+m}{2}\rfloor} \frac{(-1)^k(4^k-1)\zeta(2k+1)}{(p+m-2k)!(2\pi)^{2k}}\bigg)\\ -\sum_{k=1}^{\lfloor \frac{m-1}{2} \rfloor} \frac{(-1)^k\zeta(2k+1)}{(m-1-2k)!(p+m-2k)\pi^{2k}}-\frac{1}{(m-1)!(p+m)^2},
\end{multline}

\vspace{0.3cm}

and

\begin{multline}
\sum_{n=1}^{\infty} \frac{\zeta(2n)}{n(2n+1)\dots(2n+m-1)(2n+m+p)16^n} = \frac{\log(\pi/2)-H_{m-1}}{(m-1)!(p+m)}+ \frac{\pi}{2}*\\\sum_{k=\lfloor \frac{m+2}{2} \rfloor}^{\lfloor\frac{p+m+1}{2}\rfloor} \frac{p!(-1)^{m+k}4^k\beta(2k)}{(p+m+1-2k)!\pi^{2k}}+\delta_{\lfloor \frac{p+m}{2} \rfloor, \frac{p+m}{2}}\frac{(-1)^{m+1}p!(-1)^{\frac{p+m}{2}}2^{p+m}}{\pi^{p+m}}\zeta(p+m+1)-\frac{p!(-1)^m}{2}*\\\sum_{k=\lfloor \frac{m+1}{2} \rfloor}^{\lfloor\frac{p+m}{2}\rfloor} \frac{(-1)^k(4^k-1)\zeta(2k+1)}{(p+m-2k)!(2\pi)^{2k}} -\sum_{k=1}^{\lfloor \frac{m-1}{2} \rfloor} \frac{(-1)^k4^k\zeta(2k+1)}{(m-1-2k)!(p+m-2k)\pi^{2k}}-\frac{1}{(m-1)!(p+m)^2}.
\end{multline}

\vspace{0.5cm}

\textbf{Remark.} When $m=1$ and $p=0$, one has

\begin{equation*}
\sum_{n=1}^{\infty} \frac{\zeta(2n)}{n(2n+1)4^n} = \log\pi-1,
\end{equation*}

\vspace{0.3cm}

and

\begin{equation*}
\sum_{n=1}^{\infty} \frac{\zeta(2n)}{n(2n+1)16^n} = \frac{2G}{\pi}-1+\log\Big(\frac{\pi}{2}\Big),
\end{equation*}

\vspace{0.4cm}

the first of which is a famous series (see \cite{TylerChernoff}) and both are immediate from $(3)$ with $\theta=\pi$ and $\theta=\pi/2$, respectively. Below we compute other sums for certain $m$ and $p$.

$$ \sum_{n=1}^{\infty} \frac{\zeta(2n)}{n(2n+3)4^n} = -\frac{3\zeta(3)}{2\pi^2}+\frac{1}{3}\log\pi-\frac{1}{9}$$
$$ \sum_{n=1}^{\infty} \frac{\zeta(2n)}{n(2n+1)(n+1)4^n} = \frac{7\zeta(3)}{2\pi^2}-\log\pi-\frac{3}{2} $$
$$ \sum_{n=1}^{\infty} \frac{\zeta(2n)}{n(2n+1)(n+1)(2n+3)4^n} = \frac{2\zeta(3)}{\pi^2}+\frac{1}{3}\log\pi-\frac{11}{18} $$
$$ \sum_{n=1}^{\infty} \frac{\zeta(2n)}{n(2n+1)(n+1)(n+2)4^n} = \frac{2\zeta(3)}{\pi^2}+\frac{31\zeta(5)}{4\pi^4}+\frac{1}{2}\log\pi-\frac{7}{8} $$
$$ \sum_{n=1}^{\infty} \frac{\zeta(2n)}{n(2n+1)\dots(2n+5)4^n} = \frac{\zeta(3)}{6\pi^2}-\frac{\zeta(5)}{\pi^4}+\frac{1}{120}\log\pi-\frac{137}{7200}$$
$$ \sum_{n=1}^{\infty} \frac{\zeta(2n)}{n(n+1)16^n} = -\frac{35\zeta(3)}{4\pi^2}+\frac{4G}{\pi}+\log\Big(\frac{\pi}{2}\Big)-\frac{1}{2} $$
$$ \sum_{n=1}^{\infty} \frac{\zeta(2n)}{n(2n+1)(n+1)16^n} = \frac{35\zeta(3)}{4\pi^2}+\log\Big(\frac{\pi}{2}\Big)-\frac{3}{2} $$
$$ \sum_{n=1}^{\infty} \frac{\zeta(2n)}{n(2n+1)(2n+3)16^n} = \frac{3\zeta(3)}{8\pi^2}+\frac{8\beta(4)}{\pi^3}+\frac{1}{3}\log\Big(\frac{\pi}{2}\Big)-\frac{4}{9}$$
$$ \sum_{n=1}^{\infty} \frac{\zeta(2n)}{n(2n+1)\dots(2n+4)16^n} = \frac{2\zeta(3)}{\pi^2}-\frac{527\zeta(5)}{32\pi^4}+\frac{1}{24}\log\Big(\frac{\pi}{2}\Big)-\frac{25}{288} $$

\vspace{0.5cm}

Now, we will revisit the cotangent function and investigate a general zeta series using its power series.  

\bigskip

\section{General $\zeta(2n)$ series using $\cot(x)$}

\vspace{0.5cm}

Similar to above, we will investigate the double integral

$$ \int_{0}^{\pi z} \int_{0}^{x} x^p(x-t)^mt\cot(t) \hspace{3pt} dt \hspace{3pt} dx. $$

\vspace{0.4cm}

Using the binomial theorem, $(13)$, and a change of variables among other things,

$$ \int_{0}^{\pi z} \int_{0}^{x} x^p(x-t)^mt\cot(t) \hspace{3pt} dt \hspace{3pt} dx = \sum_{j=0}^{m} (-1)^j \binom{m}{j}\int_{0}^{\pi z} x^{p+m-j} \int_{0}^{x} t^{j+1}\cot(t) \hspace{3pt} dt \hspace{3pt} dx $$
\begin{multline*} 
= \sum_{j=0}^{m} \sum_{k=0}^{j+1} \binom{m}{j}\frac{(-1)^j(j+1)!(-1)^{\lfloor \frac{k+3}{2} \rfloor}}{(j+1-k)!2^k}\int_{0}^{\pi z} x^{p+m+1-k}\Clausen_{k+1}(2x) \hspace{3pt} dx \\ + \sum_{j=0}^{m} \delta_{\lfloor \frac{j+1}{2} \rfloor, \frac{j+1}{2}}\frac{(-1)^j(j+1)!(-1)^{\frac{j+1}{2}}}{2^{j+1}}\binom{m}{j}\zeta(j+2)\int_{0}^{\pi z} x^{p+m-j} \hspace{3pt} dx
\end{multline*}
\begin{multline*}
= \sum_{j=0}^{m} \binom{m}{j} \frac{(-1)^{j+1}}{2^{m+p+2}} \int_{0}^{2\pi z} u^{p+m+1}\Clausen_{1}(u) \hspace{3pt} du + \sum_{k=0}^{m}\sum_{j=k}^{m} \binom{m}{j} \frac{(-1)^j(j+1)!(-1)^{\lfloor \frac{k}{2} \rfloor}}{(j-k)!2^{p+m+2}}*\\\int_{0}^{2\pi z} u^{p+m-k}\Clausen_{k+2}(u) \hspace{3pt} du - \sum_{k=1}^{\lfloor \frac{m+1}{2} \rfloor} \frac{(2k)(-1)^km!(\pi z)^{p+m+2-2k}}{2^{2k}(m+1-2k)!(p+m+2-2k)}\zeta(2k+1)
\end{multline*}
\begin{multline*}
=-\frac{\delta_{m,0}}{2^{p+2}}\int_{0}^{2\pi z} u^{p+1}\Clausen_{1}(u) \hspace{3pt} du - (1-\delta_{m,0})\frac{(-1)^{\lfloor \frac{m}{2} \rfloor}m!}{2^{p+m+2}}\int_{0}^{2\pi z} u^{p+1}\Clausen_{m+1}(u) \hspace{3pt} du \\ + \frac{(m+1)!(-1)^{\lfloor \frac{m+1}{2} \rfloor}}{2^{p+m+2}}\int_{0}^{2\pi z} u^p\Clausen_{m+2}(u) \hspace{3pt} du + \sum_{k=1}^{\lfloor \frac{m+1}{2} \rfloor} \frac{(2k)(-1)^{k+1}m!(\pi z)^{p+m+2-2k}\zeta(2k+1)}{2^{2k}(m+1-2k)!(p+m+2-2k)}.
\end{multline*}

\vspace{0.5cm}

where we've evaluated the sum on $j$ and $k$. The terms with $\delta_{m,0}$ cancel each other. We can use $(21)$ to evaluate the other integrals. Rearranging and simplifying, we will arrive at

\begin{multline}
\int_{0}^{\pi z} \int_{0}^{x} x^p(x-t)^mt\cot(t) \hspace{3pt} dt \hspace{3pt} dx = (p+m+2)\bigg(2(-1)^m(\pi z)^{m+p+3}p!m!*\\\sum_{k=m+3}^{p+m+3}\frac{(-1)^{\lfloor \frac{k}{2} \rfloor}\Clausen_{k}(2\pi z)}{(p+m+3-k)!(2\pi z)^k}+\delta_{\lfloor \frac{p+m}{2} \rfloor, \frac{p+m}{2}}\frac{p!(-1)^{\frac{p+m}{2}}(-1)^mm!\zeta(p+m+3)}{2^{m+p+2}}\bigg)\\-\frac{(-1)^{\lfloor \frac{m+1}{2} \rfloor}m!(\pi z)^{p+1}\Clausen_{m+2}(2\pi z)}{2^{m+1}}+\sum_{k=1}^{\lfloor \frac{m+1}{2} \rfloor} \frac{(2k)(-1)^{k+1}m!(\pi z)^{p+m+2-2k}\zeta(2k+1)}{2^{2k}(m+1-2k)!(p+m+2-2k)}.
\end{multline}

\vspace{0.5cm}

Another way to evaluate this double integral is to use $(1)$ and Fubini's theorem, similar to the first section. Doing so, we see

\begin{multline*} 
\int_{0}^{\pi z} \int_{0}^{x} x^p(x-t)^mt\cot(t) \hspace{3pt} dt \hspace{3pt} dx = -2\int_{0}^{\pi z} x^p\sum_{n=0}^{\infty} \frac{\zeta(2n)}{\pi^{2n}}\int_{0}^{x} (x-t)^mt^{2n} \hspace{3pt} dt \hspace{3pt} dx\\= -2\sum_{n=0}^{\infty} \frac{\zeta(2n)m!\Gamma(2n+1)}{\pi^{2n}\Gamma(m+2n+2)}\int_{0}^{\pi z} x^{2n+m+p+1} \hspace{3pt} dx\\=-2\sum_{n=0}^{\infty} \frac{\zeta(2n)m!(\pi z)^{m+p+2}z^{2n}}{(2n+1)\dots(2n+m+1)(2n+m+p+2)}. 
\end{multline*}

\vspace{0.5cm}

Putting this and $(31)$ together, we find

\begin{multline}
\sum_{n=0}^{\infty} \frac{\zeta(2n)z^{2n}}{(2n+1)\dots(2n+m+1)(2n+m+p+2)} = (p+m+2)\bigg((-1)^{m+1}p!*\\\sum_{k=m+3}^{p+m+3}\frac{\pi z(-1)^{\lfloor \frac{k}{2} \rfloor}\Clausen_{k}(2\pi z)}{(p+m+3-k)!(2\pi z)^k}-\delta_{\lfloor \frac{p+m}{2} \rfloor, \frac{p+m}{2}}\frac{p!(-1)^{\frac{p+m}{2}}(-1)^m\zeta(p+m+3)}{2(2\pi z)^{m+p+2}}\bigg)\\+\frac{(-1)^{\lfloor \frac{m+1}{2} \rfloor}\Clausen_{m+2}(2\pi z)}{2(2\pi z)^{m+1}}+\sum_{k=1}^{\lfloor \frac{m+1}{2} \rfloor} \frac{k(-1)^k\zeta(2k+1)}{(2\pi z)^{2k}(m+1-2k)!(m+p+2-2k)}.
\end{multline}

\vspace{0.4cm}

When $p=0$, this formula simplifies to

\begin{multline}
\sum_{n=0}^{\infty} \frac{\zeta(2n)z^{2n}}{(2n+1)\dots(2n+m+2)} = \frac{(-1)^{\lfloor \frac{m+1}{2} \rfloor}\Clausen_{m+2}(2\pi z)}{2(2\pi z)^{m+1}}\\+\frac{(-1)^{\lfloor \frac{m}{2} \rfloor}\big(\Clausen_{m+3}(2\pi z)-\delta_{\lfloor \frac{m}{2} \rfloor, \frac{m}{2}}\zeta(m+3)\big)(m+2)}{2(2\pi z)^{m+2}}+\sum_{k=1}^{\lfloor \frac{m+1}{2} \rfloor} \frac{k(-1)^k\zeta(2k+1)}{(2\pi z)^{2k}(m+2-2k)!}.
\end{multline}

\vspace{0.4cm}

Setting $z=1/2$ and $z=1/4$, we obtain

\begin{multline}
\sum_{n=0}^{\infty} \frac{\zeta(2n)}{(2n+1)\dots(2n+m+2)4^n} = -\delta_{\lfloor \frac{m+1}{2} \rfloor, \frac{m+1}{2}}\frac{(-1)^{\frac{m+1}{2}}(2^{m+1}-1)\zeta(m+2)}{2(2\pi)^{m+1}}\\ -\delta_{\lfloor \frac{m}{2} \rfloor, \frac{m}{2}}\frac{(-1)^{\frac{m}{2}}(2^{m+3}-1)(m+2)\zeta(m+3)}{2(2\pi)^{m+2}}+\sum_{k=1}^{\lfloor \frac{m+1}{2} \rfloor} \frac{k(-1)^k\zeta(2k+1)}{(2\pi z)^{2k}(m+2-2k)!},
\end{multline}

\vspace{0.3cm}

and 

\begin{multline}
\sum_{n=0}^{\infty} \frac{\zeta(2n)}{(2n+1)\dots(2n+m+2)16^n} = (-1)^{\lfloor \frac{m}{2} \rfloor}(m+2)*\\\bigg(\delta_{\lfloor \frac{m+1}{2} \rfloor, \frac{m+1}{2}}\frac{2^{m+1}\beta(m+2)}{\pi^{m+2}}-\delta_{\lfloor \frac{m}{2} \rfloor, \frac{m}{2}}\frac{(2^{2m+5}+2^{m+2}-1)\zeta(m+3)}{4(2\pi)^{m+2}}\bigg)+(-1)^{\lfloor \frac{m+1}{2} \rfloor}*\\\bigg(\delta_{\lfloor \frac{m}{2} \rfloor, \frac{m}{2}}\frac{2^m\beta(m+2)}{\pi^{m+1}}-\delta_{\lfloor \frac{m+1}{2} \rfloor, \frac{m+1}{2}}\frac{(2^{m+1}-1)\zeta(m+2)}{4(2\pi)^{m+1}}\bigg)+\sum_{k=1}^{\lfloor \frac{m+1}{2} \rfloor} \frac{k(-4)^k\zeta(2k+1)}{\pi^{2k}(m+2-2k)!}.
\end{multline}

\vspace{0.4cm}

For general $p$,

\begin{multline}
\sum_{n=0}^{\infty} \frac{\zeta(2n)}{(2n+1)\dots(2n+m+1)(2n+m+p+2)4^n} = \frac{(-1)^mp!(m+p+2)}{2}*\\\Bigg(\sum_{k=\lfloor \frac{m+3}{2} \rfloor}^{\lfloor \frac{p+m+2}{2} \rfloor}\frac{(-1)^k(4^k-1)\zeta(2k+1)}{(p+m+2-2k)!(2\pi)^{2k}}-\delta_{\lfloor \frac{p+m}{2} \rfloor, \frac{p+m}{2}}\frac{(-1)^{\frac{p+m}{2}}\zeta(p+m+3)}{\pi^{m+p+2}}\Bigg)\\-\delta_{\lfloor \frac{m+1}{2} \rfloor, \frac{m+1}{2}}\frac{(-1)^{\frac{m+1}{2}}(2^{m+1}-1)\zeta(m+2)}{2(2\pi)^{m+1}}+\sum_{k=1}^{\lfloor \frac{m+1}{2} \rfloor} \frac{k(-1)^k\zeta(2k+1)}{\pi^{2k}(m+1-2k)!(m+p+2-2k)},
\end{multline}

\vspace{0.4cm}

and

\begin{multline}
\sum_{n=0}^{\infty} \frac{\zeta(2n)}{(2n+1)\dots(2n+m+1)(2n+m+p+2)16^n} =\frac{(-1)^mp!(m+p+2)}{4}*\\\Bigg(\sum_{k=\lfloor \frac{m+3}{2} \rfloor}^{\lfloor \frac{p+m+2}{2} \rfloor}\frac{(-1)^k(4^k-1)\zeta(2k+1)}{(p+m+2-2k)!(2\pi)^{2k}}-\sum_{k=\lfloor \frac{m+4}{2} \rfloor}^{\lfloor \frac{p+m+3}{2} \rfloor} \frac{(-1)^k\beta(2k)4^k}{(p+m+3-2k)!\pi^{2k-1}}\\-\delta_{\lfloor \frac{p+m}{2} \rfloor, \frac{p+m}{2}}\frac{(-1)^{\frac{p+m}{2}}\zeta(p+m+3)2^{p+m+3}}{\pi^{m+p+2}}\Bigg)+\delta_{\lfloor \frac{m}{2} \rfloor, \frac{m}{2}} \frac{(-1)^{\lfloor \frac{m+1}{2} \rfloor}2^m\beta(m+2)}{\pi^{m+1}}\\-\delta_{\lfloor \frac{m+1}{2} \rfloor, \frac{m+1}{2}}\frac{(-1)^{\frac{m+1}{2}}(2^{m+1}-1)\zeta(m+2)}{4(2\pi)^{m+1}}+\sum_{k=1}^{\lfloor \frac{m+1}{2} \rfloor} \frac{k(-1)^k\zeta(2k+1)4^k}{\pi^{2k}(m+1-2k)!(m+p+2-2k)}.
\end{multline}

\vspace{0.5cm}

\textbf{Remark.}
For $m=0$ and $p=0$, we have

$$\sum_{n=0}^{\infty} \frac{\zeta(2n)}{(2n+1)(n+1)4^n} = -\frac{7\zeta(3)}{2\pi^2},$$

\vspace{0.3cm}

and

$$\sum_{n=0}^{\infty} \frac{\zeta(2n)}{(2n+1)(n+1)16^n} = \frac{2G}{\pi}-\frac{35\zeta(3)}{4\pi^2},$$

\vspace{0.4cm}

the first of which was rediscovered by Ewell (see \cite{Ewell}). Below we compute other sums for certain $m$ and $p$.

$$\sum_{n=0}^{\infty} \frac{\zeta(2n)}{(2n+1)(2n+3)4^n} = -\frac{9\zeta(3)}{8\pi^2}$$
$$\sum_{n=0}^{\infty} \frac{\zeta(2n)}{(2n+1)(n+2)4^n} = -\frac{3\zeta(3)}{\pi^2}+\frac{31\zeta(5)}{2\pi^4}$$
$$ \sum_{n=0}^{\infty} \frac{\zeta(2n)}{(2n+1)(2n+2)(2n+3)4^n} = -\frac{5\zeta(3)}{8\pi^2} $$
$$ \sum_{n=0}^{\infty} \frac{\zeta(2n)}{(2n+1)(n+1)(n+2)4^n} = -\frac{\zeta(3)}{2\pi^2}-\frac{31\zeta(5)}{2\pi^4} $$
$$ \sum_{n=0}^{\infty} \frac{\zeta(2n)}{(2n+1)\dots(2n+5)4^n} = -\frac{\zeta(3)}{6\pi^2}+\frac{49\zeta(5)}{32\pi^4} $$
$$\sum_{n=0}^{\infty} \frac{\zeta(2n)}{(2n+1)(2n+3)16^n} = -\frac{9\zeta(3)}{16\pi^2}+\frac{G}{\pi}-\frac{12\beta(4)}{\pi^3}$$
$$ \sum_{n=0}^{\infty} \frac{\zeta(2n)}{(2n+1)(2n+2)(2n+3)16^n} = -\frac{61\zeta(3)}{16\pi^2} + \frac{12\beta(4)}{\pi^3} $$
$$ \sum_{n=0}^{\infty} \frac{\zeta(2n)}{(2n+1)\dots(2n+4)16^n} = -\frac{2\zeta(3)}{\pi^2}+\frac{527\zeta(5)}{16\pi^4}-\frac{4\beta(4)}{\pi^3} $$

\vspace{0.4cm}

Now, we will use these results to establish the same families of general rational zeta series but for $\zeta(2n+1)$. 

\vspace{0.5cm}

\section{Rational $\zeta(2n+1)$ series using $\psi(x)$}

\vspace{0.3cm}

\bt{For $p \in \mathbb{N}$ and $|z|<1$,}

\begin{equation}
\int_{0}^{z} x^p\psi(x) \hspace{3pt} dx = \sum_{k=0}^{p} \frac{p!(-1)^k\psi^{(-k-1)}(z)}{(p-k)!}z^{p-k}.
\end{equation}
\et

\vspace{0.3cm}

\textit{Proof.} Let $f(z)$ be the left hand side and $g(z)$ be the right hand side. It is clear that $f(0)=g(0)=0$ and $f'(z)=z^p\psi(z)$. Taking the derivative of $g(z)$, we find

$$ g'(z) = \sum_{k=0}^{p} \frac{p!(-1)^k\psi^{(-k)}(z)}{(p-k)!}z^{p-k} + \sum_{k=0}^{p-1} \frac{p!(-1)^k\psi^{(-k-1)}(z)}{(p-k-1)!}z^{p-k-1} $$
$$ = z^p\psi(z) + \sum_{k=1}^{p} \frac{p!(-1)^k\psi^{(-k)}(z)}{(p-k)!}z^{p-k} + \sum_{k=0}^{p-1} \frac{p!(-1)^k\psi^{(-k-1)}(z)}{(p-k-1)!}z^{p-k-1}. $$

\vspace{0.4cm}

Reindexing the first sum, one can see the two sums cancel each other out. So, $g'(z) = f'(z)$ for all $z$. Since they are equal at $z=0$, then $f(z) = g(z)$. $\square$

\bigskip

We can also compute this integral using $(12)$. So we will have

$$ \int_{0}^{z} x^p\psi(x) \hspace{3pt} dx = \int_{0}^{z} x^p\bigg(-\gamma-\frac{1}{x}+\sum_{k=2}^{\infty} (-1)^k\zeta(k)x^{k-1}\bigg) \hspace{3pt} dx$$
$$ = -\frac{\gamma}{p+1}z^{p+1}-\frac{z^p}{p}+\sum_{k=2}^{\infty} \frac{(-1)^k\zeta(k)}{k+p}z^{k+p}.$$

\vspace{0.4cm}

Setting these two results equal to each other, one has

\begin{equation}
\sum_{k=2}^{\infty} \frac{(-1)^k\zeta(k)}{k+p}z^k = \frac{1}{p}+\frac{\gamma z}{p+1} + \sum_{k=0}^{p} \frac{p!(-1)^k\psi^{(-k-1)}(z)}{(p-k)!z^k}.
\end{equation}

\vspace{0.5cm}

Now we can split this up as follows:

$$ \sum_{k=2}^{\infty} \frac{(-1)^k\zeta(k)z^k}{k+p} = \sum_{n=1}^{\infty} \frac{\zeta(2n)z^{2n}}{2n+p} - \sum_{n=1}^{\infty} \frac{\zeta(2n+1)z^{2n+1}}{2n+p+1}. $$

\vspace{0.4cm}

Using $(13)$ and rearranging, we see

\begin{multline}
\sum_{n=1}^{\infty} \frac{\zeta(2n+1)z^{2n}}{2n+p+1}
= -\frac{1}{2zp} - \frac{\gamma}{p+1} +\sum_{k=0}^{p} \frac{p!(-1)^{\lfloor \frac{k+1}{2} \rfloor}\pi}{(p-k)!(2\pi z)^{k+1}}\Clausen_{k+1}(2\pi z) \\- \delta_{\lfloor \frac{p}{2} \rfloor, \frac{p}{2}}\frac{p!(-1)^{\frac{p}{2}}\pi}{(2\pi z)^{p+1}}\zeta(p+1)- \sum_{k=0}^{p} \frac{p!(-1)^k\psi^{(-k-1)}(z)}{(p-k)!z^{k+1}}.
\end{multline}

\vspace{0.5cm}

For $z=1/2$ and $z=1/4$, we find

\begin{multline}
\sum_{n=1}^{\infty} \frac{\zeta(2n+1)}{(2n+p+1)4^n} = -\frac{1}{p} -\log2-\frac{\gamma}{p+1}-\sum_{k=1}^{\lfloor\frac{p}{2}\rfloor} \frac{p!(-1)^k(4^k-1)\zeta(2k+1)}{(p-2k)!(2\pi)^{2k}}\\-\delta_{\lfloor \frac{p}{2} \rfloor, \frac{p}{2}}\frac{p!(-1)^{\frac{p}{2}}\zeta(p+1)}{\pi^p}-2\sum_{k=0}^{p} \frac{p!(-2)^k\psi^{(-k-1)}(1/2)}{(p-k)!},
\end{multline}

\vspace{0.4cm}

and

\begin{multline}
\sum_{n=1}^{\infty} \frac{\zeta(2n+1)}{(2n+p+1)16^n} = -\frac{2}{p}-\log2-\frac{\gamma}{p+1}-\sum_{k=1}^{\lfloor\frac{p}{2}\rfloor} \frac{p!(-1)^k(4^k-1)\zeta(2k+1)}{(p-2k)!(2\pi)^{2k}} \\+ \pi\sum_{k=1}^{\lfloor\frac{p+1}{2}\rfloor} \frac{p!(-4)^k\beta(2k)}{(p+1-2k)!\pi^{2k}}-\delta_{\lfloor \frac{p}{2} \rfloor, \frac{p}{2}}\frac{p!(-1)^{\frac{p}{2}}2^{p+1}\zeta(p+1)}{\pi^p}-4\sum_{k=0}^{p} \frac{p!(-4)^k\psi^{(-k-1)}(1/4)}{(p-k)!}.
\end{multline}

\vspace{0.5cm}

Below we compute a few examples for specific $p$. Note that $A$ is the Glaisher-Kinkelin constant, defined by $\log A = \frac{1}{12}-\zeta'(-1)$.

$$ \sum_{n=1}^{\infty} \frac{\zeta(2n+1)}{(n+1)4^n} = -2-\gamma+12\log A - \frac{1}{3}\log2 $$
$$ \sum_{n=1}^{\infty} \frac{\zeta(2n+1)}{(2n+3)4^n} = -\frac{1}{2}-\frac{\gamma}{3}+4\log A - \frac{1}{3}\log2 $$
$$ \sum_{n=1}^{\infty} \frac{\zeta(2n+1)}{(2n+5)4^n} = -\frac{199}{180}-\frac{\gamma}{5}+8\log A -56\zeta'(-3)-\frac{3}{5}\log2 $$
$$ \sum_{n=1}^{\infty} \frac{\zeta(2n+1)}{(2n+2)16^n} = -2 - \frac{\gamma}{2} +18\log A + \log2\pi^2 - 4\log\Gamma\Big(\frac{1}{4}\Big) $$
$$ \sum_{n=1}^{\infty} \frac{\zeta(2n+1)}{(2n+3)16^n} = -64\zeta'\Big(-2,\frac{1}{4}\Big)+\frac{1}{2}+4\log A + \log2\pi^2-\frac{\gamma}{3}+\frac{3\zeta(3)}{2\pi^2}-4\log\Gamma\Big(\frac{1}{4}\Big) $$

\vspace{0.4cm}

\section{General $\zeta(2n+1)$ series using $\psi^{(-m)}(x)$}

\vspace{0.3cm}

\bt{For $p \in \mathbb{N}_{0}$, $m \in \mathbb{N}$ and $|z|<1$,}

\begin{equation}
\int_{0}^{z} x^p\psi^{(-m)}(x) \hspace{3pt} dx = \sum_{k=0}^{p} \frac{p!(-1)^k\psi^{(-k-m-1)}(z)}{(p-k)!}z^{p-k}.
\end{equation}
\et

\vspace{0.3cm}

\textit{Proof.} The proof is exactly the same as the proof of $(38)$ as there was no dependence on $\psi(x) = \psi^{(0)}(x)$. $\square$

\bigskip

Now, let us compute this integral using $(11)$ and $(12)$. Doing so, we see

$$\int_{0}^{z} x^p\psi^{(-m)}(x) \hspace{3pt} dx = \int_{0}^{z} \frac{x^p}{(m-2)!}\int_{0}^{x} (x-t)^{m-2}\bigg(-\gamma t -\log t +\sum_{k=2}^{\infty} \frac{(-1)^k\zeta(k)}{k}t^k\bigg) \hspace{3pt} dt \hspace{3pt} dx$$
$$ = \int_{0}^{z} \frac{-x^p}{(m-2)!}\bigg(\frac{\gamma x^m}{m(m-1)}+\frac{x^{m-1}(\log x-H_{m-1})}{m-1}-\sum_{k=2}^{\infty} \frac{(-1)^k\zeta(k)\Gamma(m-1)x^{m+k-1}}{k(k+1)\dots(k+m-1)}\bigg) \hspace{3pt} dx$$
\begin{multline*}
= -\frac{\gamma z^{m+p+1}}{m!(m+p+1)}+\frac{z^{m+p}(H_{m-1}-\log z)}{(m-1)!(m+p)}+\frac{z^{m+p}}{(m-1)!(m+p)^2}\\+\sum_{k=2}^{\infty} \frac{(-1)^k\zeta(k)z^{m+p+k}}{k(k+1)\dots(k+m-1)(k+m+p)}.
\end{multline*}

\vspace{0.5cm}

Setting this result equal to $(43)$ and simplifying, one has the nice result

\begin{multline}
\sum_{k=2}^{\infty} \frac{(-1)^k\zeta(k)z^k}{k(k+1)\dots(k+m-1)(k+m+p)} = \frac{\gamma z}{m!(m+p+1)}\\+\frac{\log z - H_{m-1}}{(m-1)!(m+p)}-\frac{1}{(m-1)!(m+p)^2}+\sum_{k=0}^{p} \frac{p!(-1)^k\psi^{(-k-m-1)}(z)}{(p-k)!z^{m+k}}.
\end{multline}

\vspace{0.4cm}

In the special case of $p=0$, we have

\begin{equation}
\sum_{k=2}^{\infty} \frac{(-1)^k\zeta(k)z^k}{k(k+1)\dots(k+m)} = \frac{\gamma z}{(m+1)!}+\frac{\log z - H_{m}}{m!}+\frac{\psi^{(-m-1)}(z)}{z^m}.
\end{equation}

\vspace{0.4cm}

Again, we can split $(44)$ into

\begin{multline*}
\sum_{k=2}^{\infty} \frac{(-1)^k\zeta(k)z^k}{k(k+1)\dots(k+m-1)(k+m+p)} = \sum_{n=1}^{\infty} \frac{\zeta(2n)z^{2n}}{2n(2n+1)\dots(2n+m-1)(2n+m+p)}\\ - \sum_{n=1}^{\infty} \frac{\zeta(2n+1)z^{2n+1}}{(2n+1)\dots(2n+m)(2n+1+m+p)},
\end{multline*}

\vspace{0.4cm}

and the first sum has been computed earlier. Using $(25)$, we find

\begin{multline}
\sum_{n=1}^{\infty} \frac{\zeta(2n+1)z^{2n}}{(2n+1)\dots(2n+m)(2n+m+1+p)} = -\sum_{k=0}^{p} \frac{p!(-1)^k\psi^{(-k-m-1)}(z)}{(p-k)!z^{m+k+1}}\\+\frac{(-1)^mp!}{2z}\bigg(\sum_{k=m}^{m+p} \frac{(-1)^{\lfloor \frac{k+1}{2} \rfloor}}{(p+m-k)!(2\pi z)^k}\Clausen_{k+1}(2\pi z) - \delta_{\lfloor \frac{p+m}{2} \rfloor, \frac{p+m}{2}}\frac{(-1)^{\frac{p+m}{2}}\zeta(p+m+1)}{(2\pi z)^{p+m}}\bigg)\\-\frac{1}{2z}\sum_{k=1}^{\lfloor \frac{m-1}{2} \rfloor} \frac{(-1)^k\zeta(2k+1)}{(m-1-2k)!(m+p-2k)(2\pi z)^{2k}}+\frac{\log(2\pi/z)+H_{m-1}}{2z(m-1)!(p+m)}\\+\frac{1}{2z(m-1)!(p+m)^2}-\frac{\gamma}{m!(p+m+1)}.
\end{multline}

\vspace{0.4cm}

If $p=0$, we have the nice representation

\begin{multline}
\sum_{n=1}^{\infty} \frac{\zeta(2n+1)z^{2n}}{(2n+1)\dots(2n+m+1)} = -\frac{\psi^{(-m-1)}(z)}{z^{m+1}}-\frac{1}{2z}\sum_{k=1}^{\lfloor \frac{m-1}{2} \rfloor} \frac{(-1)^k\zeta(2k+1)}{(m-2k)!(2\pi z)^{2k}}\\+\frac{(-1)^{\lfloor \frac{m}{2} \rfloor}(\Clausen_{m+1}(2\pi z) - \delta_{\lfloor \frac{m}{2} \rfloor, \frac{m}{2}}\zeta(m+1))}{2z(2\pi z)^m}+\frac{\log(2\pi/z)+H_{m}}{2zm!}-\frac{\gamma}{(m+1)!}.
\end{multline}

\vspace{0.3cm}

and for $z=1/2$ and $z=1/4$, we have

\begin{multline}
\sum_{n=1}^{\infty} \frac{\zeta(2n+1)}{(2n+1)\dots(2n+m+1)4^n} = \frac{\log4\pi +H_{m}}{m!} -\frac{\gamma}{(m+1)!}\\-\delta_{\lfloor \frac{m}{2} \rfloor, \frac{m}{2}}\frac{(-1)^{\frac{m}{2}}(2^{m+1}-1)\zeta(m+1)}{(2\pi)^m} -\sum_{k=1}^{\lfloor \frac{m-1}{2} \rfloor} \frac{(-1)^k\zeta(2k+1)}{(m-2k)!\pi^{2k}}-2^{m+1}\psi^{(-m-1)}(1/2),
\end{multline}

\vspace{0.3cm}

and

\begin{multline}
\sum_{n=1}^{\infty} \frac{\zeta(2n+1)}{(2n+1)\dots(2n+m+1)16^n} =\frac{2\log8\pi +2H_{m}}{m!} -\frac{\gamma}{(m+1)!} \\-\delta_{\lfloor \frac{m+1}{2} \rfloor,\frac{m+1}{2}} \frac{(-1)^{\frac{m+1}{2}}2^{m+1}\beta(m+1)}{\pi^m} - \delta_{\lfloor \frac{m}{2} \rfloor, \frac{m}{2}} \frac{(-1)^{\frac{m}{2}}(2^{2m+1}+2^m-1)\zeta(m+1)}{(2\pi)^m}\\-\sum_{k=1}^{\lfloor \frac{m-1}{2} \rfloor} \frac{2(-4)^k\zeta(2k+1)}{(m-2k)!\pi^{2k}}-4^{m+1}\psi^{(-m-1)}(1/4).
\end{multline}

\vspace{0.4cm}

For general $p$,

\begin{multline}
\sum_{n=1}^{\infty} \frac{\zeta(2n+1)}{(2n+1)\dots(2n+m)(2n+1+m+p)4^n} = \frac{\log4\pi +H_{m-1}}{(m-1)!(m+p)} \\+\frac{1}{(m-1)!(m+p)^2}-\frac{\gamma}{m!(m+p+1)}-(-1)^mp!\sum_{k=\lfloor \frac{m+1}{2} \rfloor}^{\lfloor\frac{p+m}{2}\rfloor} \frac{(-1)^k(4^k-1)\zeta(2k+1)}{(p+m-2k)!(2\pi)^{2k}} \\+\delta_{\lfloor \frac{p+m}{2} \rfloor, \frac{p+m}{2}}\frac{(-1)^{m+1}p!(-1)^{\frac{p+m}{2}}}{\pi^{p+m}}\zeta(p+m+1)-2^{m+1}\sum_{k=0}^{p} \frac{p!(-2)^k\psi^{(-k-m-1)}(1/2)}{(p-k)!}\\-\sum_{k=1}^{\lfloor \frac{m-1}{2} \rfloor} \frac{(-1)^k\zeta(2k+1)}{(m-1-2k)!(m+p-2k)\pi^{2k}},
\end{multline}

\vspace{0.3cm}

and

\begin{multline}
\sum_{n=1}^{\infty} \frac{\zeta(2n+1)}{(2n+1)\dots(2n+m)(2n+1+m+p)16^n} = \frac{2}{(m-1)!}\bigg(\frac{\log8\pi +H_{m-1}}{(m+p)} \\+\frac{1}{(m+p)^2}\bigg)+\sum_{k=\lfloor \frac{m+2}{2} \rfloor}^{\lfloor\frac{p+m+1}{2}\rfloor} \frac{p!(-1)^{m+k}4^k\beta(2k)}{(p+m+1-2k)!\pi^{2k-1}} -\sum_{k=\lfloor \frac{m+1}{2} \rfloor}^{\lfloor\frac{p+m}{2}\rfloor} \frac{p!(-1)^{k+m}(4^k-1)\zeta(2k+1)}{(p+m-2k)!(2\pi)^{2k}}\\-\delta_{\lfloor \frac{p+m}{2} \rfloor, \frac{p+m}{2}}\frac{p!(-1)^m(-1)^{\frac{p+m}{2}}2^{p+m+1}\zeta(p+m+1)}{\pi^{p+m}}-\frac{\gamma}{m!(m+p+1)} \\-\sum_{k=1}^{\lfloor \frac{m-1}{2} \rfloor} \frac{2(-4)^k\zeta(2k+1)}{(m-1-2k)!(m+p-2k)\pi^{2k}}-4^{m+1}\sum_{k=0}^{p} \frac{p!(-4)^k\psi^{(-k-m-1)}(1/4)}{(p-k)!}.
\end{multline}

\vspace{0.4cm}

Below we compute some sums for certain $m$ and $p$.

$$ \sum_{n=1}^{\infty} \frac{\zeta(2n+1)}{(2n+1)(n+1)4^n} = -12\log A +2 -\gamma+\frac{7}{3}\log2$$
$$ \sum_{n=1}^{\infty} \frac{\zeta(2n+1)}{(2n+1)(2n+3)4^n} = -2\log A +\frac{1}{4}-\frac{\gamma}{3}+\frac{2}{3}\log2$$
$$ \sum_{n=1}^{\infty} \frac{\zeta(2n+1)}{(2n+1)(n+2)4^n} = -4\log A +20\zeta'(-3)+\frac{19}{36}-\frac{\gamma}{2}-\frac{89}{90}\log2$$
$$ \sum_{n=1}^{\infty} \frac{\zeta(2n+1)}{(2n+1)(n+1)(2n+3)4^n} = -8\log A +\frac{3}{2}-\frac{\gamma}{3}+\log2$$
$$ \sum_{n=1}^{\infty} \frac{\zeta(2n+1)}{(2n+1)(n+1)(n+2)4^n} = -8\log A -20\zeta'(-3)+\frac{53}{36}-\frac{\gamma}{2}+\frac{121}{90}\log2 $$
$$\sum_{n=1}^{\infty} \frac{\zeta(2n+1)}{(2n+1)\dots(2n+5)4^n} = -\frac{2}{3}\log A +\frac{8}{3}\zeta'(-3)+\frac{551}{4320}-\frac{\gamma}{120}+\frac{1}{24}\log2$$
$$\sum_{n=1}^{\infty} \frac{\zeta(2n+1)}{(2n+1)(n+1)16^n} = -36\log A + 4 -\gamma+8\log2$$
$$\sum_{n=1}^{\infty} \frac{\zeta(2n+1)}{(2n+1)(2n+3)16^n} = 32\zeta'\Big(-2,\frac{1}{4}\Big)-2\log A-\frac{3}{4\pi^2}\zeta(3)-\frac{1}{4}-\frac{\gamma}{3}+2\log2$$
$$\sum_{n=1}^{\infty} \frac{\zeta(2n+1)}{(2n+1)\dots(2n+4)16^n} = -8\log A +45\zeta'(-3)+\frac{187}{144}-\frac{\gamma}{24}+\frac{41}{60}\log2$$

\vspace{0.5cm}

To conclude, we will revisit the digamma function and find another general $\zeta(2n+1)$ series.

\section{General $\zeta(2n+1)$ series using $\psi(x)$}

\vspace{0.4cm}

Similar to section 4, we will investigate the double integral

$$ \int_{0}^{z} \int_{0}^{x} x^p(x-t)^mt\psi(t) \hspace{3pt} dt \hspace{3pt} dx. $$

\vspace{0.3cm}

Using the binomial theorem and $(38)$ among other things, we find

$$
\int_{0}^{z} \int_{0}^{x} x^p(x-t)^mt\psi(t) \hspace{3pt} dt \hspace{3pt} dx = \sum_{j=0}^{m}(-1)^j\binom{m}{j}\int_{0}^{z} x^{p+m-j} \int_{0}^{x} t^{j+1}\psi(t) \hspace{3pt} dt \hspace{3pt} dx
$$

$$
= \sum_{j=0}^{m} \sum_{k=0}^{j+1}\binom{m}{j}\frac{(-1)^j(j+1)!(-1)^k}{(j+1-k)!}\int_{0}^{z} x^{p+m+1-k}\psi^{(-k-1)}(x) \hspace{3pt} dx
$$

\begin{multline*}
= \sum_{j=0}^{m}\binom{m}{j}(-1)^j\int_{0}^{z} x^{p+m+1} \psi^{(-1)}(x) \hspace{3pt} dx  \\+ \sum_{k=0}^{m} \sum_{j=k}^{m} \binom{m}{j}\frac{(-1)^j(j+1)!(-1)^{k+1}}{(j-k)!}\int_{0}^{z}x^{p+m-k}\psi^{(-k-2)}(x) \hspace{3pt} dx
\end{multline*}

\begin{multline*}
=\delta_{m,0}\int_{0}^{z} x^{p+1}\psi^{(-1)}(x) \hspace{3pt} dx + (1-\delta_{m,0})m!\int_{0}^{z} x^{p+1}\psi^{(-m-1)}(x) \hspace{3pt} dx - (m+1)!\int_{0}^{z} x^p\psi^{(-m-2)}(x) \hspace{3pt} dx
\end{multline*}

\vspace{0.4cm}

Simplifying and using $(43)$, we find

\begin{equation} 
\int_{0}^{z} \int_{0}^{x} x^p(x-t)^mt\psi(t) \hspace{3pt} dt \hspace{3pt} dx = m!z^p \bigg(z\psi^{(-m-2)}(z) - (m+p+2)\sum_{k=0}^{p} \frac{p!(-1)^k\psi^{(-k-m-3)}(z)}{(p-k)!z^k}\bigg).
\end{equation}

\vspace{0.4cm}

Now we will evaluate the same integral using the power series for $\psi(x)$ and Fubini's theorem once again. Doing so, we see

$$ \int_{0}^{z} \int_{0}^{x} x^p(x-t)^mt\psi(t) \hspace{3pt} dt \hspace{3pt} dx = -\int_{0}^{z} x^p\int_{0}^{x} (x-t)^m\bigg(t\gamma+1-\sum_{k=2}^{\infty} (-1)^k\zeta(k)t^k\bigg) \hspace{3pt} dt \hspace{3pt} dx $$
$$ = -\int_{0}^{z} x^p\Bigg(\frac{x^{m+2}\gamma}{(m+1)(m+2)}+\frac{x^{m+1}}{m+1}-\sum_{k=2}^{\infty} \frac{(-1)^k\zeta(k)m!x^{k+m+1}}{(k+1)\dots(k+m+1)}\Bigg) \hspace{3pt} dx $$
\begin{multline*}
= \frac{-z^{p+m+3}\gamma}{(m+1)(m+2)(m+p+3)} + \frac{-z^{p+m+2}}{(m+1)(m+p+2)}\\+\sum_{k=2}^{\infty} \frac{(-1)^k\zeta(k)m!z^{p+m+2+k}}{(k+1)\dots(k+m+1)(k+m+p+2)}
\end{multline*}

\vspace{0.4cm}

Using this result and $(52)$, we arrive at

\begin{multline}
\sum_{k=2}^{\infty} \frac{(-1)^k\zeta(k)z^k}{(k+1)\dots(k+m+1)(k+m+p+2)} = \frac{\gamma z}{(m+2)!(p+m+3)}\\+\frac{1}{(m+1)!(p+m+2)}+\frac{\psi^{(-m-2)}(z)}{z^{m+1}}-(m+p+2)\sum_{k=0}^{p} \frac{p!(-1)^k\psi^{(-k-m-3)}(z)}{(p-k)!z^{k+m+2}}.
\end{multline}

\vspace{0.4cm}

Note when $p=0$, 

\begin{equation}
\sum_{k=2}^{\infty} \frac{(-1)^k\zeta(k)z^k}{(k+1)\dots(k+m+2)} = \frac{\gamma z+m+3}{(m+3)!}+\frac{\psi^{(-m-2)}(z)}{z^{m+1}}-\frac{(m+2)\psi^{(-m-3)}(z)}{z^{m+2}}.
\end{equation}

\vspace{0.4cm}

As before, we will split the sum in $(53)$ as

\begin{multline*}
\sum_{k=2}^{\infty} \frac{(-1)^k\zeta(k)z^k}{(k+1)\dots(k+m+1)(k+m+p+2)} = \sum_{n=1}^{\infty} \frac{\zeta(2n)z^{2n}}{(2n+1)\dots(2n+m+1)(2n+m+p+2)} \\- \sum_{n=1}^{\infty} \frac{\zeta(2n+1)z^{2n+1}}{(2n+2)\dots(2n+m+2)(2n+m+p+3)}
\end{multline*}

\vspace{0.4cm}

and using $(32)$, we have

\begin{multline}
\sum_{n=1}^{\infty} \frac{\zeta(2n+1)z^{2n}}{(2n+2)\dots(2n+m+2)(2n+m+p+3)}= (p+m+2)\bigg((-1)^{m+1}p!*\\\sum_{k=m+3}^{p+m+3}\frac{\pi(-1)^{\lfloor \frac{k}{2} \rfloor}\Clausen_{k}(2\pi z)}{(p+m+3-k)!(2\pi z)^k}-\delta_{\lfloor \frac{p+m}{2} \rfloor, \frac{p+m}{2}}\frac{p!(-1)^{\frac{p+m}{2}}(-1)^m\zeta(p+m+3)}{2z(2\pi z)^{m+p+2}}\\+\sum_{k=0}^{p} \frac{p!(-1)^k\psi^{(-k-m-3)}(z)}{(p-k)!z^{k+m+3}}\bigg)-\frac{\gamma}{(m+2)!(p+m+3)}-\frac{1}{2z(m+1)!(p+m+2)}\\+\frac{(-1)^{\lfloor \frac{m+1}{2} \rfloor}\Clausen_{m+2}(2\pi z)}{2z(2\pi z)^{m+1}}-\frac{\psi^{(-m-2)}(z)}{z^{m+2}}+\sum_{k=1}^{\lfloor \frac{m+1}{2} \rfloor} \frac{k(-1)^k\zeta(2k+1)}{z(2\pi z)^{2k}(m+1-2k)!(m+p+2-2k)}.
\end{multline}

\vspace{0.4cm}

If $p=0$, this simplifies to

\begin{multline}
\sum_{n=1}^{\infty} \frac{\zeta(2n+1)z^{2n}}{(2n+2)\dots(2n+m+3)} = \frac{(-1)^{\lfloor \frac{m+1}{2} \rfloor}\Clausen_{m+2}(2\pi z)}{2z(2\pi z)^{m+1}}\\+\frac{(-1)^{\lfloor \frac{m}{2} \rfloor}\big(\Clausen_{m+3}(2\pi z)-\delta_{\lfloor \frac{m}{2} \rfloor, \frac{m}{2}}\zeta(m+3)\big)(m+2)}{2z(2\pi z)^{m+2}}+\sum_{k=1}^{\lfloor \frac{m+1}{2} \rfloor} \frac{k(-1)^k\zeta(2k+1)}{z(2\pi z)^{2k}(m+2-2k)!}\\-\frac{(2\gamma z+m+3)}{2z(m+3)!}-\frac{z\psi^{(-m-2)}(z)-(m+2)\psi^{(-m-3)}(z)}{z^{m+3}}.
\end{multline}

\vspace{0.4cm}

For $z=1/2$ and $z=1/4$, we see

\begin{multline}
\sum_{n=1}^{\infty} \frac{\zeta(2n+1)}{(2n+2)\dots(2n+m+3)4^n} = -\delta_{\lfloor \frac{m+1}{2} \rfloor, \frac{m+1}{2}}\frac{(-1)^{\frac{m+1}{2}}(2^{m+1}-1)\zeta(m+2)}{(2\pi)^{m+1}}\\ -\delta_{\lfloor \frac{m}{2} \rfloor, \frac{m}{2}}\frac{(-1)^{\frac{m}{2}}(2^{m+3}-1)(m+2)\zeta(m+3)}{(2\pi)^{m+2}}+\sum_{k=1}^{\lfloor \frac{m+1}{2} \rfloor} \frac{2k(-1)^k\zeta(2k+1)}{\pi^{2k}(m+2-2k)!}\\-2^{m+2}\psi^{(-m-2)}(1/2)+(m+2)2^{m+3}\psi^{(-m-3)}(1/2)-\frac{\gamma}{(m+3)!}-\frac{1}{(m+2)!},
\end{multline}

\vspace{0.3cm}

and

\begin{multline}
\sum_{n=1}^{\infty} \frac{\zeta(2n+1)}{(2n+2)\dots(2n+m+3)16^n} = (m+2)\bigg(4^{m+3}\psi^{(-m-3)}(1/4)+(-1)^{\lfloor \frac{m}{2} \rfloor}*\\\bigg(\delta_{\lfloor \frac{m+1}{2} \rfloor, \frac{m+1}{2}}\frac{2^{m+3}\beta(m+3)}{\pi^{m+2}}-\delta_{\lfloor \frac{m}{2} \rfloor, \frac{m}{2}}\frac{(2^{2m+5}+2^{m+2}-1)\zeta(m+3)}{(2\pi)^{m+2}}\bigg)\bigg)+(-1)^{\lfloor \frac{m+1}{2} \rfloor}*\\\bigg(\delta_{\lfloor \frac{m}{2} \rfloor, \frac{m}{2}}\frac{2^{m+2}\beta(m+2)}{\pi^{m+1}}-\delta_{\lfloor \frac{m+1}{2} \rfloor, \frac{m+1}{2}}\frac{(2^{m+1}-1)\zeta(m+2)}{(2\pi)^{m+1}}\bigg)-4^{m+2}\psi^{(-m-2)}(1/4)\\+4\sum_{k=1}^{\lfloor \frac{m+1}{2} \rfloor} \frac{k(-4)^k\zeta(2k+1)}{\pi^{2k}(m+2-2k)!}-\frac{\gamma}{(m+3)!}-\frac{2}{(m+2)!}.
\end{multline}

\vspace{0.3cm}

For general $p$,

\begin{multline}
\sum_{n=1}^{\infty} \frac{\zeta(2n+1)}{(2n+2)\dots(2n+m+2)(2n+m+p+3)4^n} = (-1)^mp!(m+p+2)*\\\Bigg(\sum_{k=\lfloor \frac{m+3}{2} \rfloor}^{\lfloor \frac{p+m+2}{2} \rfloor}\frac{(-1)^k(4^k-1)\zeta(2k+1)}{(p+m+2-2k)!(2\pi)^{2k}}-\delta_{\lfloor \frac{p+m}{2} \rfloor, \frac{p+m}{2}}\frac{(-1)^{\frac{p+m}{2}}\zeta(p+m+3)}{\pi^{m+p+2}}+8(-2)^m*\\\sum_{k=0}^{p} \frac{(-2)^k\psi^{(-k-m-3)}(1/2)}{(p-k)!}\Bigg)-\delta_{\lfloor \frac{m+1}{2} \rfloor, \frac{m+1}{2}}\frac{(-1)^{\frac{m+1}{2}}(2^{m+1}-1)\zeta(m+2)}{(2\pi)^{m+1}}-2^{m+2}\psi^{(-m-2)}(1/2)\\+\sum_{k=1}^{\lfloor \frac{m+1}{2} \rfloor} \frac{2k(-1)^k\zeta(2k+1)}{\pi^{2k}(m+1-2k)!(m+p+2-2k)}-\frac{\gamma}{(m+2)!(p+m+3)}-\frac{1}{(m+1)!(p+m+2)},
\end{multline}

\vspace{0.4cm}

and

\begin{multline}
\sum_{n=1}^{\infty} \frac{\zeta(2n+1)}{(2n+2)\dots(2n+m+2)(2n+m+p+3)16^n} =(-1)^mp!(m+p+2)*\\\Bigg(\sum_{k=\lfloor \frac{m+3}{2} \rfloor}^{\lfloor \frac{p+m+2}{2} \rfloor}\frac{(-1)^k(4^k-1)\zeta(2k+1)}{(p+m+2-2k)!(2\pi)^{2k}}-\sum_{k=\lfloor \frac{m+4}{2} \rfloor}^{\lfloor \frac{p+m+3}{2} \rfloor} \frac{(-1)^k\beta(2k)4^k}{(p+m+3-2k)!\pi^{2k-1}}+64(-4)^m*\\\sum_{k=0}^{p} \frac{(-4)^k\psi^{(-k-m-3)}(1/4)}{(p-k)!}-\delta_{\lfloor \frac{p+m}{2} \rfloor, \frac{p+m}{2}}\frac{(-1)^{\frac{p+m}{2}}\zeta(p+m+3)2^{p+m+3}}{\pi^{m+p+2}}\Bigg)-4^{m+2}\psi^{(-m-2)}(1/4)\\+\delta_{\lfloor \frac{m}{2} \rfloor, \frac{m}{2}} \frac{(-1)^{\lfloor \frac{m+1}{2} \rfloor}2^{m+2}\beta(m+2)}{\pi^{m+1}}-\frac{\gamma}{(m+2)!(p+m+3)}-\frac{2}{(m+1)!(p+m+2)}\\-\delta_{\lfloor \frac{m+1}{2} \rfloor, \frac{m+1}{2}}\frac{(-1)^{\frac{m+1}{2}}(2^{m+1}-1)\zeta(m+2)}{(2\pi)^{m+1}}+\sum_{k=1}^{\lfloor \frac{m+1}{2} \rfloor} \frac{k(-1)^k\zeta(2k+1)4^{k+1}}{\pi^{2k}(m+1-2k)!(m+p+2-2k)}.
\end{multline}

\vspace{0.5cm}

Below we compute some sums for certain $m$ and $p$.

$$ \sum_{n=1}^{\infty} \frac{\zeta(2n+1)}{(n+1)(2n+3)4^n} = 4\log A - 1 -\frac{\gamma}{3}+\frac{1}{3}\log2$$

$$ \sum_{n=1}^{\infty} \frac{\zeta(2n+1)}{(n+1)(n+2)4^n} = 60\zeta'(-3)-\frac{5}{12}-\frac{\gamma}{2}+\frac{19}{30}\log2$$

$$ \sum_{n=1}^{\infty} \frac{\zeta(2n+1)}{(n+1)(2n+5)4^n} = -\frac{4}{3}\log A +\frac{112\zeta'(-3)}{3} + \frac{19}{270} - \frac{\gamma}{5} + \frac{13}{45}\log2 $$

$$ \sum_{n=1}^{\infty} \frac{\zeta(2n+1)}{(n+1)(2n+3)(n+2)4^n} = 8\log A - 60\zeta'(-3)-\frac{19}{12}-\frac{\gamma}{6}+\frac{1}{30}\log2 $$

$$ \sum_{n=1}^{\infty} \frac{\zeta(2n+1)}{(2n+2)(2n+3)(2n+5)4^n} = \frac{4}{3}\log A - \frac{28\zeta'(-3)}{3}-\frac{289}{1080}-\frac{\gamma}{30}+\frac{1}{90}\log2 $$

$$ \sum_{n=1}^{\infty} \frac{\zeta(2n+1)}{(2n+2)\dots(2n+5)4^n} = \frac{2}{3}\log A - \frac{17\zeta'(-3)}{3}-\frac{277}{2160}-\frac{\gamma}{120}-\frac{1}{360}\log2 $$

$$ \sum_{n=1}^{\infty} \frac{\zeta(2n+1)}{(n+1)(2n+3)16^n} = 128\zeta'\Big(-2,\frac{1}{4}\Big)+28\log A-\frac{3\zeta(3)}{\pi^2}-5-\frac{\gamma}{3} $$

$$
\sum_{n=1}^{\infty} \frac{\zeta(2n+1)}{(n+1)(n+2)16^n} = 384\zeta'\Big(-2,\frac{1}{4}\Big)+24\log A+540\zeta'(-3)-\frac{9\zeta(3)}{\pi^2}-\frac{41}{12}-\frac{\gamma}{2}+\frac{1}{5}\log2$$

$$ \sum_{n=1}^{\infty} \frac{\zeta(2n+1)}{(2n+2)\dots(2n+4)16^n} = -32\zeta'\Big(-2,\frac{1}{4}\Big)+8\log A - 135\zeta'(-3)+\frac{3\zeta(3)}{4\pi^2}-\frac{79}{48}-\frac{\gamma}{24}-\frac{1}{20}\log2 $$

\bigskip

\begin{flushright}
\begin{minipage}{148mm}\sc\footnotesize
University of Pittsburgh, Department of Mathematics, 301 Thackeray Hall, Pittsburgh, PA 15260, USA\\
{\it E--mail address}: {\tt djo15@pitt.edu} \vspace*{3mm}
\end{minipage}
\end{flushright}

\end{document}